\documentclass[letterpaper, 10 pt, conference]{ieeeconf}
\usepackage{cite}
\usepackage{amsmath,amssymb,amsfonts}
\usepackage{algorithmic}
\usepackage{graphicx}
\usepackage{textcomp}
\usepackage{xcolor}
\usepackage{framed}
\usepackage{ulem}
\usepackage{soul}
\usepackage{changes}
\usepackage{amsmath}
\usepackage{mdframed}
\usepackage{hyperref}
\usepackage{epstopdf}

\newcommand{\calH}{{\mathcal{H}}}

\newcommand{\calU}{{\mathcal{U}}}

\newcommand{\calX}{{\mathcal{X}}}
\newcommand{\calY}{{\mathcal{Y}}}
\newcommand{\calV}{{\mathcal{V}}}

\newcommand{\KK}{{\mathbb{K}}}

\newcommand{\RR}{{\mathbb{R}}}
\newcommand{\NN}{{\mathbb{N}}}

\newcommand{\knl}{\mathfrak{K}}

\newcommand{\unl}{\ell}

\newtheorem{theorem}{Theorem}

\newtheorem{lemma}{Lemma}

\def\BibTeX{{\rm B\kern-.05em{\sc i\kern-.025em b}\kern-.08em
    T\kern-.1667em\lower.7ex\hbox{E}\kern-.125emX}}

\begin{document}
\IEEEoverridecommandlockouts
	\allowdisplaybreaks
\title{Convergence Rates of Online Critic Value Function \\Approximation   in Native Spaces
}

\author{Shengyuan  Niu$^1$, Ali Bouland$^1$, Haoran Wang$^1$, Filippos Fotiadis$^2$, \\Andrew Kurdila$^1$, Andrea L'Afflitto$^3$, Sai Tej Paruchuri$^4$, Kyriakos G.  Vamvoudakis$^2$
\thanks{$^1$S. Niu, A. Bouland, H. Wang, and A. Kurdila are with the Department of Mechanical Engineering, Virginia Tech, Blacksburg, VA, USA. Email:
{\tt\small \{syniu97, bouland, haoran9, kurdila\}@vt.edu}.}
\thanks{$^2$F. Fotiadis and K. G. Vamvoudakis are with The Daniel Guggenheim School of Aerospace Engineering, 
Georgia Institute of Technology, Atlanta, GA, 
USA. Email:
{\tt\small \{ffotiadis, kyriakos\}@gatech.edu}.}
\thanks{$^3$A.  L'Afflitto is with the Grado Department of Industrial and Systems Engineering, Virginia Tech, Blacksburg, VA, USA. Email: {\tt\small a.lafflitto@vt.edu}.}
\thanks{$^{4}$S. T. Paruchuri is with the Department of Mechanical Engineering and Mechanics, Lehigh University, Bethlehem, PA, USA. Email: {\tt\small saitejp@lehigh.edu}. 
}
\thanks{
This work was supported in part, by NSF under grants
No. CAREER CPS-$1851588$,
CPS-$2227185$,
S\&AS-$1849198$, and
$2137159$, and
the US Army Research Lab under Grant No. $W911QX2320001$.}
}
\maketitle

\begin{abstract}

\textcolor{black}{This paper derives rates of convergence of online critic methods for the estimation of the value function for a class of nonlinear optimal control problems.} 
Assuming that the underlying value function lies in reproducing kernel Hilbert space (RKHS),
we derive explicit bounds on the performance of the critic in terms of the kernel functions, 
the number of basis functions, and
the scattered location of centers used to define the RKHS.
The performance of the critic is precisely measured in terms of the power function of the scattered bases, and
it can be used either in an {\it a priori} evaluation of potential bases or in an {\it a posteriori} assessments of the value function error for basis enrichment or pruning.
The most concise bounds in the paper describe explicitly how the critic performance depends on the placement of centers,
as measured by their fill distance in a subset that contains the trajectory of the critic. \textcolor{black}{To the authors' knowledge, precise error bounds of this form are the first of their kind for online critic formulations used in optimal control problems. In addition to their general and immediate applicability to a wide range of applications, they have the potential to constitute the groundwork for more advanced ``basis-adaptive'' methods for nonlinear optimal control strategies, ones that address limitations due to the dimensionality of approximations.}  
\end{abstract}


\section{Introduction}

\textcolor{black}{Optimal control has become one of the core methodologies in modern control theory for nonlinear systems. 
One of its main advantages lies in its ability to yield control laws that achieve a compromise between the control effort expended and the time needed to attain regulation.  
At the heart of nonlinear optimal control design is the Hamilton-Jacobi-Bellman (HJB) equation \cite{lewis2012optimal},
a partial differential equation (PDE) that is notoriously difficult to solve analytically. 
Numerous studies describe methods to approximate the solution of the HJB equation, see the reviews in \cite{kiumarsi2017optimal,vamvoudakis2021handbook} for example.  One of the most popular such tools is policy iteration (PI), which is a 
 process that cyclically evaluates the cost function for a given controller, and, subsequently,
improves that controller from measurements.
Nevertheless, an issue with PI is its need to employ a neural network (NN) for the policy evaluation step, called the  ``critic,''
which inherently leads to approximation errors that degrade performance or lead to failure of convergence of the PI process.}

\textcolor{black}{Notable early efforts  that study Galerkin approximations for PI in a recursive implementation include  \cite{bea1998successive,beard1997galerkin,abu2005nearly}.
Subsequent papers \cite{vamvoudakis2010online,bhasin2013novel} use some of the theory in \cite{beard1997galerkin,bea1998successive,abu2005nearly} to study various online implementations based on learning theory.  
The works referenced in \cite{lewis2013reinforcement} and \cite{kamalapurkar2018reinforcement}, for example,
provide comprehensive reviews of contemporary  theory underlying many recent online and offline methods.
Yet, these results do not derive explicit descriptions of how performance is related {\it quantitatively} to rates of convergence of value function approximations generated by a critic.
On the other hand, some very recent efforts in \cite{bian2021reinforcement,kalise2020robust,yang2021hamiltonian} emphasize the importance of examining the impact of the approximation error on the performance of reinforcement learning methods.}

\textcolor{black}{This work continues the strategy started for offline reproducing kernel Hilbert spaces (RKHS) methods in \cite{bouland2023rates}, but now considers online approaches for the critic step in PIs. 
We describe how the fill distance of the centers used to define the bases for approximation dictates the performance of the critic. In several case, we relate the rate of convergence in the RKHS directly and explicitly to the performance of the critic. To the authors' knowledge, this is the first time that such rates of convergence have been derived for the online critic step. These general results and error rates   have a host of potential applications to reduce  guesswork by the control designer when using PI techniques.  The derived rates also have the potential to serve as the foundation of methods that dynamically add or delete scattered or sparse basis functions to address cases when dimensionality of approximations becomes an issue.} 

\section{Problem Statement}
Consider the continuous-time nonlinear system
\begin{align}\label{eq:sys}
    \dot{x} (t) &= f(x(t))+g(x(t)) u (x(t)),~  x(0) = x_0, ~t\ge0,
\end{align}
where 
$x : [0,\infty) \to \mathbb{R}^n$,
$f : \mathbb{R}^n \rightarrow \mathbb{R}^n$,
$g : \mathbb{R}^n \rightarrow \mathbb{R}^{n \times m}$, and
$u : \mathbb{R}^n \rightarrow \mathbb{R}^m$ represent the
state of the system,
the drift dynamics,
the input dynamics, and
the control input, respectively.  
The problem is to find a continuous control input $t \mapsto u(t)$ that minimizes the cost functional
\begin{align}\label{eq:costJ}
    J(x_0,u) = \int_{0}^\infty \underbrace{\left( Q(x(t)) + u^\textrm{T}(t) R u(t) \right)}_{r(x(t),u(t))} \textrm{d}t
\end{align} 
where  $Q: x \mapsto  Q(x)\ge 0$, and $R \succ 0$.
One of the main issues with this problem is that one needs to solve a challenging nonlinear HJB equation.
A minimizer $u^\star$ of \eqref{eq:costJ} is called an optimal control input,
and $V^\star(\cdot)=J(\cdot,u^\star)$ defines the optimal value function.
Then, to find $u^\star$ and $V^\star$, in principle, one needs to find the positive-definite solution $V^\star$ of the HJB equation 
\begin{multline}\label{eq:HJB}
    \calH_{u^*}(V^*(x)) = 
    -\frac{1}{4}\nabla V^{\star\textrm{T}}(x)g(x)R^{-1}g^\textrm{T}(x)\nabla V^\star(x)\\+\nabla V^{\star\textrm{T}}(x)f(x)+Q(x)=0,~V^\star(0)=0,\ \forall x\in \Omega,
\end{multline}
and then
calculate $u^\star(x)=-\frac{1}{2}R^{-1}g^\textrm{T}(x)\nabla V^\star(x)$ \cite{lewis2012optimal}, where
$\calH_u(V)$ is the Hamiltonian function associated with $u$ and $V$.
Nevertheless, \eqref{eq:HJB} is generally difficult, if not impossible, to solve analytically for $V^\star$.
For this reason, PI is often employed to solve \eqref{eq:HJB} approximately \cite{abu2005nearly, bea1998successive}.

The most crucial and computationally demanding step of PI is that of policy evaluation.
Given a continuous feedback gain
$\mu:\mathbb{R}^n\rightarrow\mathbb{R}^m$ that stabilizes \eqref{eq:sys} on a set $\Omega\subseteq \mathbb{R}^n$,
policy evaluation seeks to find the value function
$V_\mu(\cdot)  \triangleq  J(\cdot,\mu)$ associated with that controller.
Provided that this function is continuously differentiable,
it follows from \cite{lewis2012optimal} that it satisfies
\begin{multline}\label{eq:LE}
    \calH_\mu(x) \triangleq  \calH_\mu(V_\mu(x))=\nabla V_\mu^\textrm{T}(x)(f(x)+g(x)\mu(x))\\+Q(x)+\mu^\textrm{T}(x)R\mu(x)=0,~V_\mu(0)=0. 
\end{multline}
While an analytical solution to \eqref{eq:LE} is also difficult to obtain,
its linearity in $V_{\mu}$ --- a property not present in \eqref{eq:HJB} --- enables the use of the so-called \textit{critic} NN as a means to approximately solve it over a compact set $\Omega\subset\mathbb{R}^n$. 

To that end, note that since $V_\mu$ is continuous,
it can be
expressed on $\Omega$ as 
$V_\mu(x)={W^\textrm{T}\phi(x)}+\epsilon_N(x),~\forall x\in\Omega$, where
$\phi:\mathbb{R}^n\rightarrow\mathbb{R}^N$ is a suitable vector of $N$ basis functions,
$W\in\mathbb{R}^N$ denote the ``ideal weights''  for that basis, and
$\epsilon_N:\mathbb{R}^n\rightarrow\mathbb{R}$ denotes the approximation error.
The critic NN then uses an estimate $\hat{W}(t)\in\mathbb{R}^N$ of $W$, and
provides an estimate of $\hat{v}_N(t,\cdot)$ of $V_\mu$ according to the formula
$\hat{v}_{N}(t,x)=\hat{W}^\textrm{T}(t)\phi(x)$.
The purpose of policy evaluation is, thus, to properly train the critic weights
$\hat{W}(t)$ so that the norm of the parameter error
$\tilde{W}(t) \triangleq  W-\hat{W}(t)$
becomes as small as possible.
In \cite{vamvoudakis2010online}, the online policy evaluation law
\begin{equation}\label{eq:critic_law}
\dot{\hat{W}}(t)=-a\frac{\sigma(t)}{(\sigma^\textrm{T}(t)\sigma(t){+}1)^2}\left(\sigma^\textrm{T}(t)\hat{W}(t){+}r(x(t),\mu(x(t))) \right)
\end{equation}
was proposed,
where
$\sigma(t) \triangleq  \sigma(x(t))=\nabla \phi(x(t))(f(x(t))+g(x(t))\mu(x(t))$, and
$a>0$ denotes the learning rate.
Interestingly, it was proved that, under a persistency of excitation condition,
the parameter estimation error $\tilde{W}(t)$ under \eqref{eq:critic_law} indeed converges exponentially fast to a neighborhood of the origin,
the size of which scales with the size of the approximation error $\epsilon_N$ over $\Omega$.
Nevertheless, the size of $\epsilon_N$ is rarely known beforehand and, to our knowledge,
no existing general strategy yet has been able to precisely quantify how the basis influences the performance of the critic. 

This paper lifts the analysis of the norm of the parameter error
$\|\hat{W}(t)-W\|_{\RR^N}$
to an analysis of 
$\|v_N(t,\cdot)-V_\mu\|_{H(\Omega)}$,
which captures 
estimates of the error of the value function.
This analysis makes explicit the contribution of approximation errors in a wide variety of choices of the RKHS $H(\Omega)$
assumed to contain the value function.
Our goal is to ultimately use this analysis to quantitatively relate the choice of the basis function $\phi$ of the critic NN to the error $\|\hat{v}_N(t,\cdot)-V_\mu\|_{H(\Omega)}$ 
in online critic estimates $v_N(t,\cdot)$ of the value function $V_\mu$. 
A further goal of the paper is to reduce trial-and-error in realistic implications of the critic for adaptive nonlinear optimal control.

 

\section{Notation and Preliminaries}
\textcolor{black}{
Denote 
$v$ as a generic value function,
$\hat{v}$ as an estimate of $v$, and
$\tilde{v}\triangleq v-\hat{v}$ as the error.} 

\subsection{Elements of RKHS Theory}
We denote by $H(\Omega)$ an RKHS over the set $\Omega\subseteq \RR^n$ that is constructed using a Mercer reproducing kernel
$\knl : \Omega \times \Omega \rightarrow \RR$.
A Mercer kernel $\knl(\cdot,\cdot)$ is continuous, symmetric, and of positive type.
Being of positive type means that, for any $N$-point subset $\Xi_N \subset \Omega$,
the corresponding Grammian matrix
$\KK_N \triangleq [\knl(\xi_i,\xi_j)]\in \RR^{N\times N}$
is positive semidefinite.
The native space $H(\Omega)$ itself is then determined as the closure of the linear span of the kernel sections $\knl_x(\cdot) \triangleq \knl(x,\cdot)$, that is,
$H(\Omega)  \triangleq  \overline{\text{span}\{\knl_x(\cdot)\ | \ x\in \Omega \}}$, 
where
the closure is taken with respect to the candidate inner product $(\knl_x,\knl_y) \triangleq  \knl(x,y)$ for all $x,y\in \Omega$. 

Approximations in this paper are constructed using the finite-dimensional subspace
$H_N \triangleq  \text{span}\{\knl_{\xi_i}(\cdot)\ | \ \xi_i\in \Xi_N, 1\leq i\leq N\}$.
We denote by $\Pi_N:H(\Omega)\to H_N$ the $H(\Omega)$-orthogonal projection of $H(\Omega)$ onto $H_N$.
\textcolor{black}{A key property of orthogonal projections onto a closed subspace of a Hilbert space is that they map an arbitrary input into the closest element of the subspace.}

The evaluation functional 
$E_x : H(\Omega) \rightarrow \RR$ is defined so that,
for each $x\in \Omega$ and every $f\in H(\Omega)$,
it holds that
$E_x f  \triangleq  f(x)$.
Thus, the evaluation functional defines the bounded linear mapping $H(\Omega) \to \RR$.
The reproducing property, which is satisfied for any RKHS,
implies that $E_xf=f(x)=(f,\knl_x)_H$ for any $f\in H(\Omega)$ and $x\in \Omega$.
Furthermore, as $E_x$ is a bounded linear operator between Hilbert spaces,
its adjoint operator $E_x^*  \triangleq  (E_x)^* : \RR \rightarrow H(\Omega)$ is a bounded linear operator.
This adjoint operator is expressed as $E_x^*\alpha  \triangleq  \knl_x \alpha$ for all $\alpha \in \RR, x\in \Omega$.
That is, $E^*_x$ can be understood as a multiplication operator since it multiplies any real number by the function $\knl_x$. 

If the kernel \( \knl(\cdot,\cdot) \) is bounded on the diagonal, then
per definition,
there exists a constant \( \bar{\knl} \) such that,
\( \knl(x, x) \leq \bar{\knl}^2 \)
for every \( x \in  \Omega \).
This condition guarantees that every function within the space \( H(\Omega) \) is continuous and bounded.
Furthermore, it ensures boundedness of  the operator norm, that is,
$\|E_x\|=\|E^*_x\|\leq \bar{\knl}$. 
It is worth noting that many commonly used kernels satisfy this criterion,
including the inverse multiquadric, Sobolev-Mat\'{e}rn, Wendland, and exponential kernels \cite{wendland}.

\subsection{Differential Operator $A$ on Native Spaces} 
\label{sec:diffop}
We begin by introducing  the differential 
 operator $A$ that is defined pointwise as 
$    (Av)(x) \triangleq  \left ( f(x) + g(x) \mu(x)  \right )^\textrm{T} \nabla v(x)$ for all $x \in \Omega,$ whenever $v$ is sufficiently smooth. 
Note that \eqref{eq:LE} then corresponds to the operator equation
$Av=b$ with
$b=-r$,
$r$ defined in terms of the kernel $r$ of the cost function in \eqref{eq:costJ}, and
$v=V_\mu$.  

\begin{theorem}
\label{th:AandAstar}
Let  the kernel $\knl:\Omega\times \Omega\to \RR$ that defines the native space $H(\Omega)$ be a $C^{2m}(\Omega,\Omega)$ function with $m\geq 1$,
and suppose that $\mu$ and $f_i,g_i$ for $1\leq i\leq d$ are multipliers for $C(\Omega)$ and $H(\Omega)$. Then,
\begin{enumerate}
    \item The  operator $A:H(\Omega)\to C(\Omega)$, as well as the operator $A:H(\Omega)\rightarrow L^2(\Omega)$,  is bounded, linear, and compact.
    \item The adjoint operator $A^*:L^2(\Omega)\rightarrow H(\Omega)$  has representation 
    \begin{align*}
        A^*&=\int_\Omega \left ( \nabla_x \knl(x,y) \right)^\mathrm{T} \left (f(x) + g(x)\mu(x) \right) h(x)\mathrm{d}x \\
& \triangleq \int_\Omega \ell^*(y,x)h(x)\mathrm{d}x
    \end{align*}
    for any $y\in \Omega$ and $h(\cdot)\in L^2(\Omega)$. 
    \item Considered as a mapping  $A^*:L^2(\Omega)\to H(\Omega)$, or as a mapping $A^*:L^2(\Omega)\rightarrow L^2(\Omega)$, the operator $A^*$ is compact. 
\end{enumerate}
\end{theorem}
\begin{proof}
The proof of this theorem can be found in \cite{bouland2023rates}, which uses  Theorem 1 of \cite{zhou2008derivative}.  \hfill  
\end{proof}

\textcolor{black}{Note that the assumptions in Theorem \ref{th:AandAstar} imply that the basis functions that define  $H_N$ are continuously differentiable.}

\subsection{The DPS Learning Law and Its Approximation}
\label{sec:online}
For developing the online learning laws, we introduce the time-varying functional 
    \[
    \mathcal{J}(t,\tilde{v}) \triangleq \frac{1}{2}|E_{x(t)}A\tilde{v}|^2=\frac{1}{2}\left (A^*E^*_{x(t)}E_{x(t)}A\tilde{v},\tilde{v} \right)_H,  
    \]
    which is defined for all $\tilde{v}\in H(\Omega)$ that satisfy the additional regularity condition 
    $
    \tilde{v}\in \{f\in H(\Omega)\ | \ A\tilde{v} \in H(\Omega)\}$. 
    The analysis in the remainder of this paper always assumes that this regularity condition holds. 
    An elementary calculation shows that the Fr\'{e}chet derivative of $\mathcal{J}(t,\tilde{v})$ is given by $D\mathcal{J} \triangleq  A^*E^*_{x(t)}E_{x(t)}A$. 
For a fixed time $t$, let $\hat{v}(t,\cdot)\in H(\Omega)$ be a time-varying approximation of the minimizer $v$   of $\mathcal{J}(t,v)$. 
An ideal gradient learning law designs the error $\tilde{v}(t,\cdot) \triangleq v-\hat{v}(t,\cdot)$ so that it evolves in the local direction of steepest descent, which is defined in terms of the Fr\'{e}chet differential in 
\begin{align*}
\frac{\partial }{\partial t} \tilde{v}(t,\cdot)=-a A^*E^*_{x(t)} ( y(t) -E_{x(t)}A\hat{v}(t,\cdot)) \in H(\Omega),
\end{align*}
where $y(t) \triangleq E_{x(t)}Av$, and $a>0$. 
This ideal gradient law evolves in $H(\Omega)$, and
defines a distributed parameter system.
In the usual way, we define the ideal evolution law for the estimate $\hat{v}(t,\cdot)$ as 
{ 
\[
\frac{\partial }{\partial t} \hat{v}(t,\cdot) 
= -a A^*E^*_{x(t)}E_{x(t)}A\hat{v}(t,\cdot) +a A^*E^*_{x(t)} y(t) \in H(\Omega).
\]
}

{\noindent Note that, in contrast to \cite{vamvoudakis2010online}, the critic state evolves in a function space and
this evolution law can be understood as a PDE.}
Finite-dimensional approximations of this PDE are obtained by choosing $\hat{v}_N(t,\cdot) \triangleq \sum_{j=1}^N \hat{W}_{j}(t) \knl_{\xi_j}(\cdot)$ and seeking a solution of 
{
\begin{align}
&\frac{\textrm{d}}{\textrm{d}t} \hat{v}_N(t,\cdot) \notag \\ 
&= -a \Pi_N A^*E^*_{x(t)}E_{x(t)}A\Pi_N \hat{v}_N(t,\cdot) +a \Pi_N A^*E^*_{x(t)} y(t). \label{eq:onlineHN}
\end{align}
}

\noindent These finite-dimensional equations evolve in $H_N$, and they are equivalent to a system of ODEs.

\subsection{Online Coordinate Realizations}
\label{sec:onlinecoords}
The critical step in deriving coordinate realizations of \eqref{eq:onlineHN} 
 must examine representations of the operator $\Pi_N A^* E_x^* E_x A \Pi_N$. The finite-dimensional approximation $\Pi_NA^*E^*_xE_xA\Pi_N$ can be deduced by considering $g=\knl_{\xi_j}$ and $h=\knl_{\xi_i}$ to obtain 
\begin{align*}
    [\mathbb{A}_N(x)]_{i,j}& \triangleq ((\Pi_NA^*E_x^*E_xA\Pi_N) \knl_{\xi_j},\knl_{\xi_i})_H ,\\ 
    &=\left[\Phi^\textrm{T}(x,\Xi_N)\psi(x)\psi(x)^\textrm{T}\Phi(x,\Xi_N)  \right ]_{i,j}.
\end{align*}
After taking the inner product of \eqref{eq:onlineHN} with an arbitrary $\knl_{\xi_i}\in H_N$, we therefore obtain the system of ODEs
\begin{align*}
    &\KK_N
    \dot{\hat{W}}(t)
    = 
    -a 
    \mathbb{A}_N(x(t))
     \hat{W}(t) + a Y(t),
\end{align*}
where $\hat{W} \triangleq  [{\hat{W}}_1(t), \ldots, {\hat{W}}_{N}(t)]^{\rm T}$,
$y(t)=E_{x(t)}Av=b(x(t))$ denotes the output,
$Y_i(t) \triangleq ( A^*E^*_{x(t)}y(t), \knl_{\xi_i})_{H(\Omega)}$ and
$Y(t)= [Y_1,\ldots, Y_N(t)]^{\rm T}$. 

  \smallskip

\noindent \underline{Remark 1}: Interestingly,
$Y(\cdot)$ is essentially equivalent to the right-hand-side of  \eqref{eq:critic_law}, with a slight difference being that the normalization with $(\sigma^\textrm{T}\sigma+1)^2$ in \eqref{eq:critic_law} is not introduced  here. 

\noindent \underline{Remark 2}: It is well-known that, in practice, the gradient learning law in \eqref{eq:onlineHN} must use a robust modification whenever external noise, numerical noise, or approximation error appears in $Y(\cdot)$.
This is the reason for the normalization ordinarily used in PI and reinforcement learning.
In the following, we discuss a dead-zone robust modification for this purpose.
\textcolor{black}{The dead-zone modification is advantageous since it enables a simpler proof of rates of convergence in some cases.}

\subsection{Rates of Convergence and Online Performance Bounds}\label{sec:onlinerates}
In our first error analysis of online algorithms, we employ the gradient learning law \eqref{eq:onlineHN}.
This analysis is based on modifying the approach in  \cite{vamvoudakis2010online} and carefully tracking the dependence of expressions on the number of bases $N$ and the approximation error.  The theorem below develops an ultimate bound on $\bar{v}_N \triangleq  \Pi_N \tilde{v}_N=\Pi_N(v-\hat{v}_N)$. 
\begin{theorem}\label{th:onlineerrorN}
Suppose that the kernel $\knl(\cdot,\cdot)$ that defines the RKHS $H(\Omega)$  is bounded on the diagonal by a constant $\bar{\knl}^2$.  In addition assume that  the family of subspaces  $\{H_N\}_{N\in \NN}$ and trajectory $t\mapsto x(t)$ are PE in the sense that there are constants $\Delta(N), \gamma_1(N)$ depending on $N$ and $\gamma_2>0$ such that 
\[
\gamma_1(N) I_{H_N} \leq \underbrace{\int_t^{t+\Delta(N)} \Pi_N A^* E^*_{x(\tau)}E_{x(\tau)}A\Pi_N \mathrm{d}\tau}_{S_N(t)} \leq \gamma_2 I_{H_N} 
\]
for each $N\in \NN$,
where $S_N(t):H_N\to H_N$.
Then, 
     \begin{align}
&\|\bar{v}_N(t,\cdot)\|_{H(\Omega)} \triangleq \|\Pi_N v-\hat{v}_N(t,\cdot)\|_{H(\Omega)} \notag\\ &\leq 
   \frac{\sqrt{\gamma_2\Delta(N)}}{\gamma_1(N)} (\bar{\calY}_{N,\max} + \delta\gamma_2 a (\bar{\calY}_{N,\max}+\epsilon_{N,\max}) ). \label{eq:VamvoudakisRate}
     \end{align}
where
\begin{align}
    \bar{\calY}_{N,\textrm{max}} &\triangleq
    \sup_{\tau\in [t,t+\Delta(N)]}|E_{x(\tau)}A\Pi_N\bar{v}_N(t,\cdot)|, \label{eq:ymax} \\
    \epsilon_{N,\textrm{max}} &\triangleq \sup_{\tau\in [t,t+\Delta(N)]} \epsilon_{N}(\tau). \label{eq:epsNmax} 
\end{align}
\end{theorem}

\begin{proof}
    The consistent approximation of the gradient law can be written as 
    \begin{align*}
    \frac{\textrm{d}}{\textrm{d}t} \bar{v}_N(t,\cdot) 
    =&- a \Pi_N A^*E^*_{x(t)} E_{x(t)} A \bar{v}_N(t,\cdot) \\
    &  - a \Pi_N A^*E^*_{x(t)} E_{x(t)} A (I-\Pi_N)v.
    \end{align*}
    Following the proof of Technical Lemma 2, part b in \cite{vamvoudakis2010online}, we rewrite this equation as the system 
    \begin{align*}
        \dot{\calX}_N(t)&=B_N(t)\calU_N(t), \\
        \calY_N(t)&=C_N^*(t) \calX_N(t),
    \end{align*}
    where $B_N(t) \triangleq -a \Pi_NA^*E^*_{x(t)}$, $C_N^*(t) \triangleq E_{x(t)}A\Pi_N$, $\calX_N(t) \triangleq \bar{v}_N(t,\cdot)$, $\epsilon_N(t) \triangleq E_{x(t)}A(I-\Pi_N)v$, and  $\calU_N(t) \triangleq -\calY_N(t)+\epsilon_N(t)$.
    Both $B_N(t)$ and $C_N(t)$ are bounded linear operators, and
    their bounds can be chosen independently of $N$ sicne the kernel $\knl(\cdot,\cdot)$ is bounded on the diagonal, and, hence,
    $\|E_{x(t)}\|=\|E^*_{x(t)}\|\leq \bar{\knl}$. 
    Also, it holds that $\calX_N(t)\in H_N$ and  $\calY_N(t)\in \RR$. 
    The proof of Technical Lemma 2 part b in \cite{vamvoudakis2010online} holds for states, controls, and observations in Euclidean spaces, like $\RR^d$ or $\RR$. 
    Since all the operators above are bounded, each step in the proof of Equation (A.9) in Technical Lemma 2 part b in \cite{vamvoudakis2010online} can also be applied without change in the current setting.
    In particular, it holds that
    \begin{align*}
    \|\calX(t)\|_{H_N} &\leq \frac{\sqrt{\gamma_2\Delta(N)}}{\gamma_1(N)}\bar{\calY}_{N,\textrm{max}}
     \\
    & \hspace*{.25in} +\frac{\delta \gamma_2}{\gamma_1(N)}
    \int_{t}^{t+\Delta(N)}\|B_N(\tau)\| \cdot \|\calU_N(\tau)\|\textrm{d}\tau
    \end{align*}
for a constant $\delta$ of order one.
    Furthermore, 
    $\|B_N(t)\|\leq a \|A^*\| \bar{\knl}$ 
     and 
     $\|\calU_N(t)\|\leq \bar{\calY}_{N,\textrm{max}} + |\epsilon_{N}(t)|$. 
    We conclude that 
 the rate in \eqref{eq:VamvoudakisRate} holds with
\begin{equation}\epsilon_{N,\textrm{max}}  \leq  \sup_{\tau\geq 0}E_{x(\tau)}A(I-\Pi_N)v.  \quad  \square 
\end{equation}
\end{proof}

The next theorem bounds the ultimate output error $\tilde{y}_N(t) \triangleq  y(t)-\hat{y}_N(t)$, where
$y(t) = E_{x(t)}Av$ and
$\hat{y}_N(t) \triangleq  E_{x(t)}A\hat{v}_N(t,\cdot)$,
in terms of the approximation error $\epsilon_{N,\textrm{max}}$ in the case whereby we use a hard dead-zone version of the learning law with a properly sized dead-zone.
We emphasize how the next result does not require a PE condition, and
the error bound on performance is more readily tied to just the approximation error $\epsilon_{N,\max}$ as described in Section \ref{sec:explicit}.
On the other hand, in principle, an oracle must define a dead-zone that is a tight bound for the approximation error.
In practice, the size of the dead-zone is defined iteratively.

\begin{theorem}\label{th:deadzone}
    Consider a learning law for 
    $\hat{v}_N(\cdot,\cdot)$, such that
    if $\tilde{y}_N(t) \triangleq  E_{x(t)}A\tilde{v}_N(t,\cdot) \geq \bar{\epsilon} \geq  \epsilon_{N,\max}$ for some $t \geq 0$, then 
    \eqref{eq:onlineHN} is verified, and,
    if $\tilde{y}_N(t) < \epsilon_{N,\max}$ for some $t \geq 0$, then
    $\frac{\textrm{d}}{\textrm{d}t} \hat{v}_N(t,\cdot) = 0$.
    Then, for any arbitrarily small constant $\eta>0$,
    there exists $T \triangleq T(\eta) > 0$ such that 
    $
    |E_{x(t)}(\calH_\mu - \hat{\calH}_N(t,\cdot))|\equiv  |\tilde{y}_{N}(t)| \leq \frac{1+\eta}{a} \bar{\epsilon} 
    $
    for all $t\geq T(\eta)$,
    where the Hamiltonian $\calH_\mu$ is defined in \eqref{eq:LE} and
    $\hat{\calH}_N(t,x) \triangleq  A\hat{v}_N(t,x)+r(x)$ denotes
    the approximate Hamiltonian with
    $x\in \Omega$.
    If we choose $\bar{\epsilon} \triangleq  M(N)\epsilon_{N,\max}$ for some (small) integer $M(N)$, and $T_O>0$ is the time that the measurement error $\tilde{y}_N(t)$ spends outside the dead-zone, then
    we an ultimate bound on the decrease of the value function error is given by 
    $\|\tilde{v}_N(t,\cdot)\|_{H(\Omega)}^2 \leq \|\tilde{v}_N(t_0,\cdot)\|_{H(\Omega)}^2 - 2aT_O(1+M(N))M(N)\epsilon_{N,\max}^2$
    for all $t\geq 0$ large enough. 
\end{theorem}
\begin{proof}
    \textcolor{black}{In this proof we choose the Lyapunov function $\calV(\tilde{v}_N) \triangleq \frac{1}{2}(\tilde{v}_N,\tilde{v}_N)_{H(\Omega)}$.
    When $|\tilde{y}_N(t)|\geq \bar{\epsilon}$,
    the derivative of the Lyapunov function along trajectories of the learning law satisfy 
    \begin{align*}
        &\frac{\textrm{d}}{\textrm{d}t}\calV(\tilde{v}_N(t,\cdot))& \\
        &= -a \left (
        E_{x(t)} A \tilde{v}_N(t,\cdot),E_{x(t)}A\tilde{v}_N(t,\cdot)\right)_{\RR}\\
        &\quad +a\left(E_{x(t)} A \tilde{v}_N(t,\cdot),-E_{x(t)}A(I-\Pi_N)\tilde{v}_N(t)\right)_{\RR}  \\
        &\leq- a |\tilde{y}_N(t)| \left ( 
        |\tilde{y}_{N}(t)|- \epsilon_{N,\max}   
        \right ).
    \end{align*}
    }
    Because $|\tilde{y}_N(t)|\geq \bar{\epsilon}\geq \epsilon_{N,\max}$, we have
    \begin{align*}
        \frac{\textrm{d}}{\textrm{d}t}\calV(\tilde{v}_N(t,\cdot))&\leq- a |\tilde{y}_N(t)| \left ( 
        |\tilde{y}_{N}(t)|- \epsilon_{N,\max}   
        \right ),\\
        &\leq - a \bar{\epsilon} \left ( 
        \bar{\epsilon}- \epsilon_{N,\max}   
        \right )<0,
    \end{align*}
while the trajectory is outside the dead-zone.
Following standard convergence arguments,
we conclude that the time spent outside the dead-zone is finite, and, thus, the norm of the output $\tilde{y}(t)$ is ultimately bounded by the dead-zone.
The bound on the value function error $\|\tilde{v}(t,\cdot)\|_{H(\Omega})$ can be derived by evaluating the Lyapunov function at the time the observations $\tilde{y}(t)$ enters the dead-zone \cite[Ch. 10]{lavretskyBook}.   \hfill   
\end{proof}

\section{Explicit Error Bounds and Fill Distances}\label{sec:explicit}
In this section, we describe how some techniques used to describe rates of convergence of approximations in a native space can be applied to the bounds in Theorems \ref{th:onlineerrorN} and \ref{th:deadzone} on the online critic.
Note that a bit more can be said about the errors $\epsilon_N(t)$ and $\epsilon_{N,\max}$ that appear in these theorems.
It holds that
\begin{align*}
    &\epsilon_N(t) {\triangleq}   |E_{x(t)}A(I{-}\Pi_N) v|  =|(\unl(\cdot,x(t)),(I-\Pi_N)v)_{H(\Omega)}|\\
    &{\leq} {\sup_{\xi\in \Omega} \|\unl(\cdot,\xi)\|_{H(\Omega)}}  \|(I{-}\Pi_N)v\|_{H(\Omega)}{\leq} \unl_{\textrm{max}} \|(I{-}\Pi_N)v\|_{H(\Omega)},
\end{align*}
where $\ell(x,y)=\ell^*(y,x)$ and $\ell^*(y,x)$ is defined in Theorem \ref{th:AandAstar}. 

The remainder of this section describes how $\|(I-\Pi_N)v\|_{H(\Omega)}$ can be explicitly bounded in terms of the placement of centers in $\Xi_N$. 
The power function of the subspace $H_N$ in the RKHS $H(\Omega)$ is defined as  
$\mathcal{P}_N(x) \triangleq  \sqrt{\knl(x,x)-\knl_N(x,x)}$
with $\knl_N$ the reproducing kernel of the subspace $H_N$
\cite{wendland,schaback94}.
It can be proven that
$\knl_N(x,y) \triangleq  \knl_{\Xi_N}(x)^\textrm{T}\KK^{-1}_N \knl_{\Xi_N}(y)$ where $\knl_\Xi(x)= [\knl_{\xi_1}(x),\ldots,\knl_{\xi_N}(x)]^\textrm{T}\in \RR^N$
denotes the column vector of $N$ basis functions defined in terms of the set of centers $\Xi_N\subset \Omega$. 
The power function is useful for generating point-wise bounds on the projection error such as
$|E_x(I-\Pi_N)v|\leq \mathcal{P}_N(x)\|(I-\Pi_N)v\|_{H(\Omega)}$
for all $x\in \Omega$ and $v\in H(\Omega)$
and any native space whatsoever \cite{wendland,schaback94}.

We use this well-known identity to bound the error $\|(I-\Pi_N)v\|_{H(\Omega)}$ that appears in the ultimate bound of the critic. 
\begin{theorem}[Modification of Theorem 11.23 in \cite{wendland}]\label{th:doubling}
    Suppose that $v$ satisfies the regularity condition
    $v = \mathcal{L}u$,
    where $\mathcal{L} : L^2(\Omega)\to H(\Omega)$ is the bounded, linear, compact operator 
        $(\mathcal{L} u)(x)  \triangleq  \int_\Omega \knl(x,y)u(y)\mathrm{d}y.$
    Then, there is a constant $C>0$ such that
    $\|(I-\Pi_N)v\|_{H(\Omega)} \leq C \sup_{\xi\in \Omega} |\mathcal{P}_N(\xi)|\|\mathcal{L}^{-1}v\|_{L^2(\Omega)}$
    provides an error bound.
\end{theorem}
\begin{proof}
    This proof is based on that of Theorem 11.23 of \cite{wendland}, and for completeness we summarize the simple modifications here. First note that 
    \begin{align*}
        (w,Lu)_{H(\Omega)}&=(w,\int_\Omega \knl(\cdot,y)u(y)\textrm{d}y)_{H(\Omega)} \\
        &=\int_\Omega (w,\knl_y)_{H(\Omega)}u(y)\textrm{d}y \\
        &= \int_\Omega w(y)u(y)\textrm{d}y=(w,u)_{L^2(\Omega)}. 
    \end{align*}
    Now we can write 
    \begin{align*}
        \|(I-\Pi_N)v\|_{H(\Omega)}^2&=((I-\Pi_N)v,(I-\Pi_N)v)_{H(\Omega)}, \\ &= ((I-\Pi_N)v,v)_{H(\Omega)}, \\
        &=((I-\Pi_N)v,Lu)_{H(\Omega)}, \\
        &=(I-\Pi_N)v,u)_{L^2(\Omega)}, \\
        &\leq \|(I-\Pi_N)v\|_{L^2(\Omega)} \|u\|_{L^2(\Omega)}.
    \end{align*}
    But we also have 
    \begin{align*}
    \|(I-\Pi_N)v\|^2_{L^2(\Omega)} &= \int_\Omega |E_x(I-\Pi_N)v|^2 \textrm{d}x \\
    &\leq |\Omega| |\sup_{\xi\in \Omega}\mathcal{P}_N(\xi)|^2\|(I-\Pi_N)v\|_{H(\Omega)}^2
    \end{align*}
    Substituting this bound above completes the proof of the theorem. 
\end{proof}
Since the centers $\Xi_N$, kernel $\knl$, and power function $\mathcal{P}_N$ are known,
Theorem \ref{th:doubling} can be used,
in either {\it a priori} or {\it a posteriori} estimation of the value function estimate error that results from using a collection of centers $\Xi$. 

The geometric nature of the bound in Theorem \ref{th:doubling}
is often emphasized by relating the power function to the fill distance of the centers $\Xi_N$ in the set $\Omega$, which
is defined as
$h_{\Xi_N,\Omega}  \triangleq  \sup_{y\in \Omega} \min_{\xi_i\in \Xi_N}\|y-\xi_i\|_2$. 
For a variety of kernel functions,
which can be applied to the problem addressed in this paper,
\cite{wendland,schaback94}
provide bounds on the power function in the form 
$\mathcal{P}_{N}(x)\lesssim \sqrt{\mathcal{N}(h_{\Xi_N,\Omega})}$, where $\mathcal{N}:\RR^+\rightarrow \RR^+$
depends on the kernel function.
The following lemma summarizes three common examples of such bounds.
\begin{lemma}\label{lem:lemma1}
Suppose that $v$ is contained in the uncertainty class 
$
\mathcal{C}_{L,R} \triangleq  \{g = \mathcal{L} u\in H(\Omega)\ | \ \|u\|_{L^2(\Omega)}\leq R\}\subseteq H(\Omega) 
$, and that the assumptions
of Theorem  \ref{th:deadzone} holds with the minimum size dead-zone $\bar{\epsilon}\approx \epsilon_{N,\max}$.
For the Sobolev-Matérn kernel of a high enough smoothness $k$ in Table \ref{tab:pf_F_bounds}, then
there exists $T>0$ such that, for all $t\geq T$,
    \begin{align*}
    &|E_{x(t)}(\calH_\mu-\hat{\calH}_N(t,\cdot))| \equiv |y(t)-\hat{y}_N(t)| \approx O(h^{k-n/2}_{\Xi_N,\Omega}). 
  \end{align*}
For the Wendland compactly supported kernel $\eta_{n,k}$,
$|E_{x(t)}(\calH_\mu-\hat{\calH}_N(t,\cdot))| \approx O(h^{k+1/2}_{\Xi_N,\Omega})$.
For the exponential kernel,
$|E_{x(t)}(\calH_\mu-\hat{\calH}_N(t,\cdot))| \approx  O\left (\sqrt{e^{-\alpha|h_{\Xi,\Omega}|/h_{\Xi_N,\Omega}}}\right)$ for a constant $\alpha$ that depends on the hyperparameters of the exponential kernel. 
\end{lemma}
\begin{proof}
    This result follows from Theorems \ref{th:deadzone} and \ref{th:doubling}. \hfill   
\end{proof}
\begin{table*}[h] 
     \centering
     \begin{tabular}{lll} \hline \hline 
     \hspace*{.5in} Type  & \hspace*{.5in} Kernel & \hspace*{.2in} $\mathcal{N}(\cdot)$ \\ 
         Gaussian & $e^{-\alpha r^2}$,$\ \ \alpha>0$ & $e^{-a|\log h|/h}$\\
          Multiquadric& $(-1)^{\lceil \beta \rceil}(c^2 + r^2)^\beta$,$\ \ \beta>0$,$\beta\not \in \NN$ & $e^{-a/h}$\\
         Inverse multiquadric & $(c^2+r^2)^\beta$, $\ \ \beta<0$ & $e^{-a/h}$\\
         Compactly supported functions& $\eta_{d,k}$ & $h^{2k+1}$ \\ 
         Sobolev-Matérn & $\frac{2\pi^{d/2}}{\Gamma(k)}{K}_{k-d/2}(r/2)^{k-d/2}$, $d,k\in \NN$ & $h^{2k-d}$ \\ 
          \hline \hline 
     \end{tabular}
     \caption{Functions $\mathcal{N}(\cdot)$ that bound the power function as  
     $\mathcal{P}_{N}(x)\lesssim \sqrt{\mathcal{N}(h_{\Xi_N,\Omega})}$ for $x\in \Omega$. In this table $h \triangleq  h_{\Xi_N,\Omega}$. The constant $a$ is a generic constant that varies with the kernel choice. The symbol $K_{\nu}$ is the modified Bessel function of the $3^{rd}$ kind of order $\nu$, while the symbol $\eta_{d,k}$ denotes the compactly supported Wendland kernel of dimension $d$ and smoothness $k$ given in \cite{wendland}.}
     \label{tab:pf_F_bounds}
 \end{table*}
\section{Numerical Results}
\label{sec:numerical}

\textcolor{black}{In this section, we present numerical validation studies for the system of the form \eqref{eq:sys} studied in  \cite{vamvoudakis2010online}, with 
$f(x)= \left[-x_1+x_2,
-0.5 x_1-0.5 x_2\left(1-\left(\cos \left(2 x_1\right)+2\right)^2\right) \right]^{\rm T}$ and $g(x) = \left[0,
\cos \left(2 x_1\right)+2 \right]^{\rm T}$.}  
The cost function for this problem sets  $R=1$ and $Q=I_2$, with $I_2$  the  identity matrix in $\RR^{2\times 2}$. The optimal  value function is
$V^\star(x) = 0.5 x_1^2 + x_2^2$, which generates  the optimal feedback controller  $u^\star(x) = -(\cos(2x_1) + 2)x_2$. 
The numerical validation studies in \cite{vamvoudakis2010online} are based on a very low-dimensional system of polynomial bases,
whose span contains the exact optimal value function.  

\begin{figure}[!hb]
    \centering
    \includegraphics[width = 1\columnwidth]{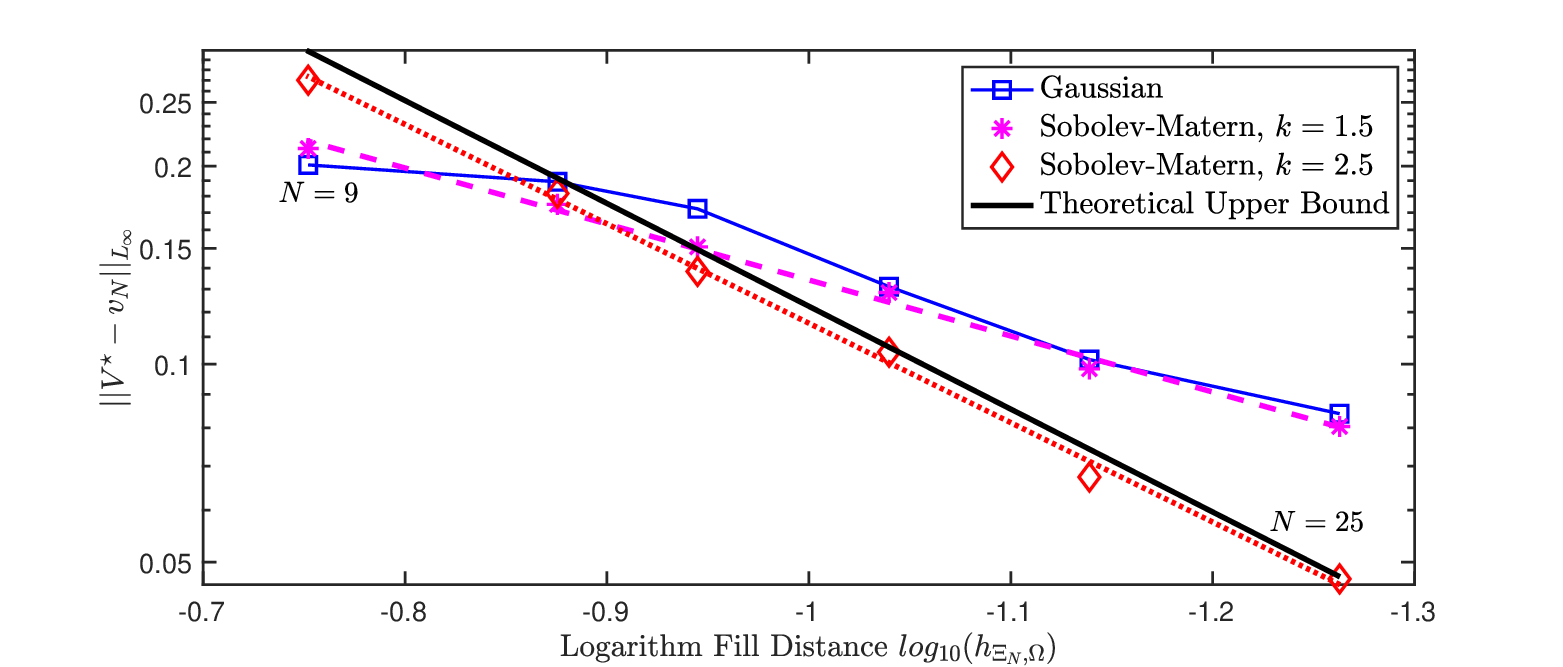}
    \caption{The $L^\infty(\Omega)$ error norm of the online critic estimates of the value function $V^\star$ using the dead-zone rule described in Lemma \ref{lem:lemma1}. The steady-state value function approximations using  Sobolev-Matérn kernels of smoothness $k=1.5, 2.5$ and exponential kernels are plotted above. Note that the rates of convergence for the Sobolev-Matérn kernels closely follow the theoretical bounds derived in Lemma \ref{lem:lemma1}.} 
    \label{fig:fig1}
\end{figure}

Figure \ref{fig:fig1} depicts the error norm $\|V^\star-\hat{v}_N(t,\cdot)\|_{L^\infty(\Omega)}$ for two  Matérn kernels and an exponential kernel.
Since $\|E_x\|\leq \bar{\knl}$,
it holds that 
$|E_x\tilde{v}_N(t,\cdot)|\leq \bar{\knl}\|\tilde{v}_N(t,\cdot)\|_{H(\Omega)}$ and $\|\tilde{v}_N(t,\cdot)\|_{L^\infty(\Omega)}\leq \bar{\knl}\|\tilde{v}_N(t,\cdot)\|_{H(\Omega)}$, and
Lemma \ref{lem:lemma1} implies the corresponding convergence in the norm of $L^\infty(\Omega)$.   
The ultimate approximate value function $\hat{v}_N(t,\cdot)$ closely matches the analytical expression for the optimal value function $V^\star$ as the dimension $N\to\infty$.
Note that Theorem \ref{lem:lemma1} only guarantees that $\hat{v}_N(t,\cdot)$ converges to $V_\mu$, not $V^\star$, and, indeed,
this plot is a more stringent empirical test of the performance of the critic. 
Figure \ref{fig:fig1} illustrates that the online critic estimates $\hat{v}_N(t,\cdot)$ for the Sobolev-Matérn kernels converge at a rate that is theoretically determined by the fill distance as described in the paper in Lemma \ref{lem:lemma1}. 

\textcolor{black}{We should emphasize that these studies make use of regular arrays of centers to verify and validate the derived error bounds.
The proposed method does not require that centers be selected using regular grids.
Suppose for example that the state trajectory of interest is embedded in a high-dimensional state space, but
its evolution resides on a low-dimensional submanifold embedded in that high-dimensional space.
Since the approach in this paper uses scattered bases that are not confined to any \textit{a priori} grids or triangulations,
their locations can be tailored to locations on the submanifold.
Such a strategy can be pursued to address issues related to the curse of dimensionality.
An in-depth study of how to execute this strategy in practice far exceeds the limits of this introductory paper.
However, in principle, the proposed method is not restricted to regular grids of bases that scale like $N^d$.}

A bit more can be deduced about the value function error in $L^\infty(\Omega)$ when the regularity condition in Lemma \ref{lem:lemma1} holds. This is referred to as the ``doubling trick'' in the literature on approximations in RKHS;
see Theorem 11.23 of \cite{wendland} that enables the conclusion $|E_x(I-\Pi_N)f|\leq  O((\sup_{\xi\in \Omega}\mathcal{P}_N(\xi))^2)$.
A line having this slope for the Sobolev-Matérn kernel with $k=2.5$ is labeled in Figure \ref{fig:fig1} as  the ``theoretical upper bound.''

Often, in implementations, it is of vital concern to establish the rates of convergence of the error $\mu-\hat{\mu}_N$ where $\hat{\mu}_N$ is the control approximation based on $\hat{v}_N(t,\cdot)$ of the ideal control $u^\star$. We can proceed exactly as in the proof of Theorem 3 of \cite{bouland2023ratesarXiv} in the case at hand to conclude that 
\begin{align*}
\|u^\star&-\hat{u}_N(t,\cdot)\|_{C(\Omega)}\leq C \|V^\star-\hat{v}_N(t,\cdot)\|_{H(\Omega)} \\
&\leq C\left ( \|V^\star-V_\mu\|_{H(\Omega)} + \|\tilde{v}_N(t,\cdot)\|_{H(\Omega} \right) 
\end{align*}
for some fixed constant $C>0$. Thus, if $\|V^*-V_\mu\|_{H(\Omega)}$ is sufficiently small, say of $O(\epsilon_{N,\max})$, then we expect the same rate of convergence for the control convergence in $C(\Omega)$ as in Lemma \ref{lem:lemma1}  for $\|\tilde{v}_N(t,\cdot)\|_{H(
\Omega)}$.

\begin{figure}[!tb]
    \centering
    \includegraphics[width = 1\columnwidth]{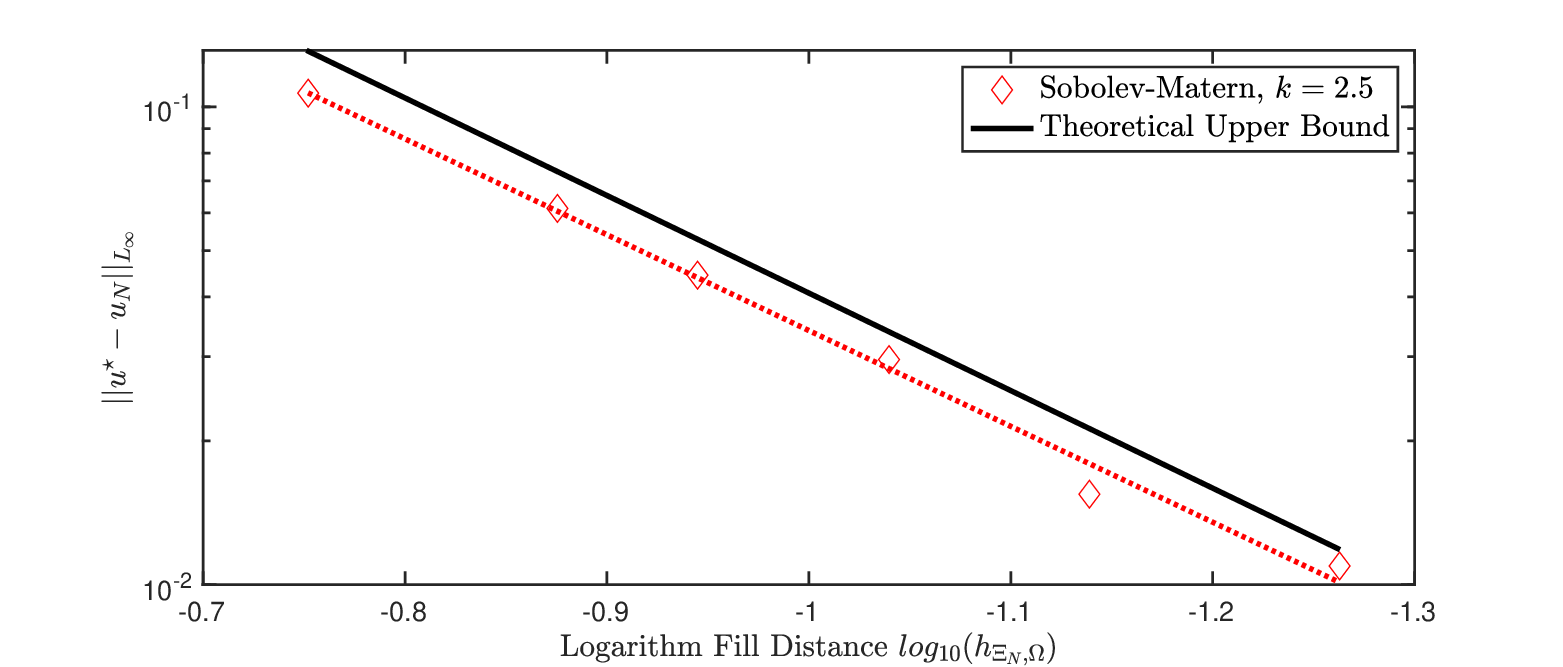}
    \caption{The $L^\infty(\Omega)$ error norm of the online critic estimates of the control input $u^\star$ with Sobolev-Matérn kernels of smoothness $k=2.5$. Note that the rates of convergence for the Sobolev-Matérn kernels closely follow the theoretical bounds derived in Theorem 3 of \cite{bouland2023ratesarXiv}.}
    \label{fig:fig2}
\end{figure}

\begin{figure}[!tb]
    \centering
    \includegraphics[width = 1\columnwidth]{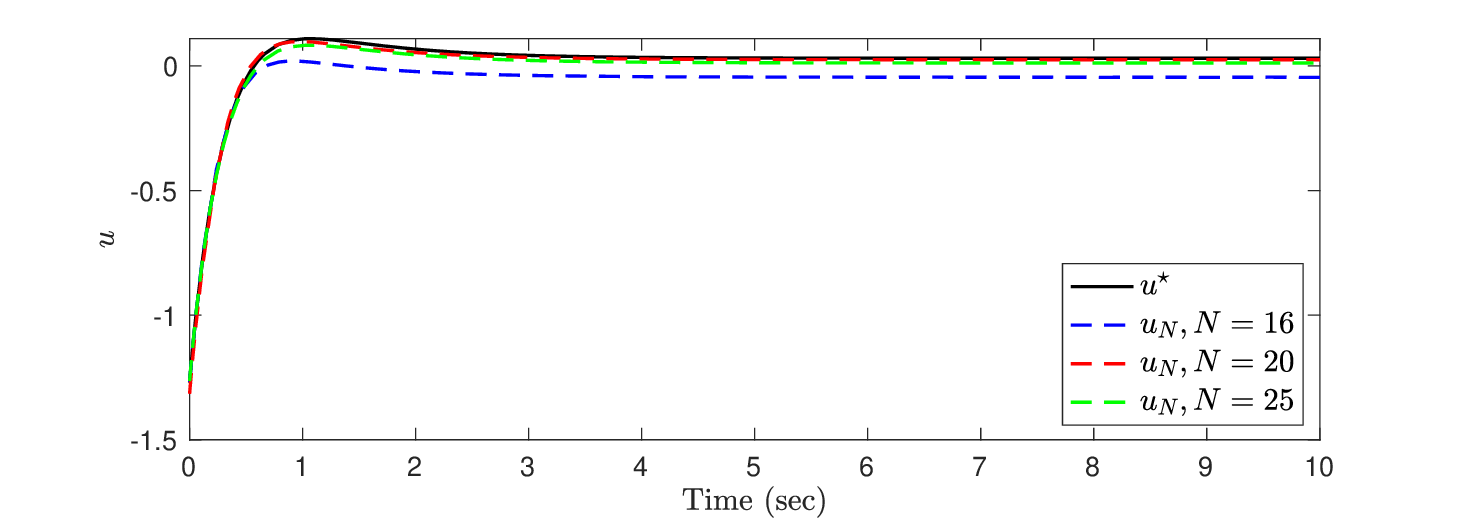}
    \caption{Feedback control $u$ by learned $\hat{W}$ with Sobolev-Matérn kernels of smoothness $k=2.5$. }
    \label{fig:fig3}
\end{figure}
\section{Conclusions}
\color{black}This paper has formulated the online critic for estimating the optimal value function in terms of evolution laws for a wide variety of RKHSs. \color{black}The paper lifts conventional approaches, \color{black}which focus on studies of the convergence of parameter errors $\|W-\hat{W}(t)\|_{\RR^N}$ in $\RR^N$, to instead focus on the norms of the value function error $\|V^\star-\hat{v}(t,\cdot)\|_{H(\Omega)}$. A wide variety  of the performance bounds on the error in the  value function estimates are derived in terms of the power function of the scattered basis. This basic result is subsequently refined to obtain performance guarantees on the critic in terms of the fill distance of the centers in the subset of interest $\Omega$.\color{black}

 \bibliographystyle{IEEEtran}
\bibliography{reference}

\end{document}